\documentclass[12pt]{amsart}

\usepackage{amssymb}
\newtheorem{thm}{Theorem}[section]
\newtheorem{theorem}[thm]{Theorem}
\newtheorem{corollary}[thm]{Corollary}
\newtheorem{lemma}[thm]{Lemma}

\theoremstyle{definition}
\newtheorem{definition}[thm]{Definition}
\newtheorem{example}[thm]{Example}

\newtheorem{remark}[thm]{Remark}
\numberwithin{equation}{section}

\newcommand{\al}{\alpha}
\newcommand{\de}{\delta}
\newcommand{\ve}{\varepsilon}
\newcommand{\ga}{\gamma}
\newcommand{\Ga}{\Gamma}
\newcommand{\ka}{\varkappa}
\newcommand{\GC}{\mathcal{G}}
\newcommand{\N}{\mathbb{N}}
\newcommand{\Om}{\Omega}
\newcommand{\PC}{\mathcal{P}}
\newcommand{\R}{\mathbb{R}}
\newcommand{\SC}{\mathcal{S}}
\newcommand{\Tf}{{\mathfrak T}}

\newcommand{\WC}{\mathcal{W}}

\newcommand{\ex}{\mathrm{ex}\,}
\newcommand{\supp}{\mathrm{supp}\,}
\newcommand{\conv}{\mathrm{conv}\,}
\newcommand{\cconv}{\overline{\mathrm{conv}}\,}

\def\leq {\leqslant}
\def\geq {\geqslant}

\begin{document}

\title[Birkhoff's theorem]
{Birkhoff's theorem\\ for a family of probability spaces}
\author{Y. Safarov}
\address{Department of Mathematics, King's College,
Strand, London, UK}
\email{ysafarov@mth.kcl.ac.uk}
\date{March 2005}
\dedicatory{To the memory of O.A. Ladyzhenskaya}


\begin{abstract}
We extend Birkhoff's theorem on doubly stochastic matrices to some
countable families of discrete probability spaces with nonempty
intersections.
\end{abstract}

\maketitle


A (possibly infinite) square matrix
$\,\{w_{ij}\}_{i,j=1,2,\ldots}\,$ with nonnegative entries
$\,w_{ij}\,$ is said to be {\it doubly stochastic} if its row and
column sums are equal to one. The matrix
$\,\{w_{ij}\}_{i,j=1,2,\ldots}\,$ can be identified with a function
$\,w\,$ on the direct product of two discrete spaces
$\,X=\{x_1,x_2,\ldots\}\,$ and $\,Y=\{y_1,y_2,\ldots\}\,$, such that
$\,w(x_i,y_j)=w_{ij}\,$, $\,\forall i,j=1,2,\ldots\,$. Under this
identification, the matrix $\,\{w_{ij}\}_{i,j=1,2,\ldots}\,$ is
doubly stochastic if and only if the restriction of $\,w\,$ to each
subset $\,\{x_i\}\times Y\,$ or $\,X\times\{y_j\}\,$ is the density
of a probability measure.

According to Birkhoff's theorem \cite{Bi1},
\begin{enumerate}
\item[(I)]
the extreme points of the convex set of doubly stochastic matrices
are permutation matrices,
\item[(II)]
the set of doubly stochastic matrices coincides with the closed
convex hull of the set of permutation matrices.
\end{enumerate}
Obviously, if (I) holds then (II) follows from
\begin{enumerate}
\item[(II$'$)]
the set of doubly stochastic matrices coincides with the closed
convex hull of the set of its extreme points.
\end{enumerate}

Many proofs of Birkhoff's theorem are known for finite matrices
(see, for example, \cite{An} \cite{BR} or \cite{Ro}). The set of
finite doubly stochastic matrices
$\,\{w_{ij}\}_{i,j=1,2,\ldots,n}\,$ is compact. Therefore (II$'$) is
a particular case of the Krein--Milman theorem.

The problem of extending (I) and (II) to infinite matrices is known
as Birkhoff's problem $111$ \cite{Bi2}. It was considered in
\cite{Is}, \cite{RP}, \cite{Ke} and \cite{Sa}. In the infinite case
(I) remains true but the validity of (II) depends on the choice of
topology. The set of doubly stochastic matrices is not compact in
any natural locally convex topology on the liner space of infinite
matrices (see Section 3). Therefore the Krein--Milman theorem is not
applicable and one has to prove (II) separately.

Birkhoff's theorem has been generalized in many directions. In
particular,
\begin{enumerate}
\item[$\bullet$]
in \cite{Ho}, \cite{Ka}, \cite{LLL}, \cite{ST} and \cite{Ti} the
authors considered various subsets of the set of doubly stochastic
matrices.
\item[$\bullet$]
The papers \cite{CLMST}, \cite{Gr}, \cite{Le} and \cite{Mu} dealt
with classes of matrices with fixed, but not necessarily equal to
one, row and column sums.
\item[$\bullet$]
A measure $\,\mu\,$ on the direct product of two unit intervals is
said to be doubly stochastic if $\,\mu(A\times X)=\mu(X\times
A)=|A|\,$ for every Borel set $\,A\subseteq X\,$, where $\,|A|\,$
denotes the Lebesgue measure of $\,A\,$. The doubly stochastic
measures were studied in \cite{BS}, \cite{Do}, \cite{Fe}, \cite{Li},
\cite{KST} and \cite{Vi}. In \cite{Do} and \cite{Li} the authors
independently obtained a continuous analogue of the first statement
of Birkhoff's theorem.
\item[$\bullet$]
In \cite{LMST} the authors considered nonnegative ``hypermatrices'',
that is, nonnegative functions defined on the direct product of
several discrete spaces.
\end{enumerate}

The aim of this paper is to obtain an analogue of Birkhoff's theorem
for a countable family of probability spaces with nonempty
intersections.  We consider a family $\,\Ga\,$ of countable sets
$\,\Om_1,\Om_2,\ldots\,$ and the convex set $\,\SC(\Ga)\,$ of
nonnegative functions defined on the union $\,\Om:=\cup_k\Om_k\,$
whose restrictions $\,\left.w\right|_{\Om_k}\,$ are densities of
probability measures on the sets $\,\Om_k\,$. Each function from
$\,\SC(\Ga)\,$ can be identified with a family of probability
measures on the sets $\,\Om_k\,$ which coincide on the intersections
$\,\Om_i\cap\Om_j\,$.

Unlike in the above mentioned papers, we do not assume that the set
$\,\Om\,$ is a direct product and the sets $\,\Om_k\,$ are its
fibres. If $\,\Om\,$ is infinite, we also consider convex subsets
$\,\SC(\Ga,\WC)\subset\SC(\Ga)\,$ which consist of functions
satisfying certain decay conditions at infinity. Such conditions
often appear in applications (see, for example, \cite{Sa}).

The idea to consider several probability spaces with nonempty
intersections is quite natural. This model seems to be a more
adequate reflection of reality than the classical scheme with one
probability space, where one implicitly assumes that all events
lying outside its scope either have probability zero or are totally
unrelated to the space under consideration.

The sets $\,\SC(\Ga,\WC)\,$ themselves and their extreme points are
determined by the layout of the sets $\,\Om_k\,$. In Section 2 we
join every two elements of $\,\Om\,$ lying in the same set
$\,\Om_k\,$ by an edge and formulate our results in terms of the
obtained graph $\,G\,$. Theorem \ref{T2.11} gives a complete
description of the set of extreme points under the assumption that
the multiplicity $\,\ka(g)\,$ of the covering
$\,\{\Om_k\}_{k=1,2,\ldots}\,$ does not exceed two for each
$\,g\in\Om\,$ (in other words, this means that each point
$\,g\in\Om\,$ belongs to at most two distinct sets $\,\Om_k\,$). If
$\,\ka(g)\leq2\,$ then each function lying in the set of extreme
points may take only the values 0, $\,\frac12\,$, 1, and the
subgraph associated with its support consists of isolated vertices
and isolated odd primitive cycles (see Definition \ref{D2.1}).

If $\,\Om\,$ is the direct product of two discrete spaces $\,X\,$
and $\,Y\,$, and $\,\Ga\,$ is the family of all sets of the form
$\,\{x_i\}\times Y\,$ or $\,X\times\{y_j\}\,$, then $\,\ka\equiv2\,$
and the corresponding graph $\,G\,$ does not have odd primitive
cycles. Therefore (I) is a special case of Theorem \ref{T2.11}. The
graphs generated by multidimensional hypermatrices contain odd
primitive cycles. Possibly, this explains the presence of the
non-standard extreme points discussed in \cite{LMST}. Even primitive
cycles were considered in \cite{Gr}, \cite{Le} and \cite{Mu}.
Theorem \ref{T2.11} shows that, under the condition
$\,\ka(g)\leq2\,$, this cycles do not affect the structure of the
set $\,\SC(\Ga)\,$.

The main result of Section 3 is Theorem \ref{T3.3}, where we prove
(II$'$) under certain conditions on topology and the family $\,G\,$.
These conditions do not include any assumptions about extreme
points. Moreover, in Section 3 we do not use any results or
definitions from Section 2. Therefore this section can be read
separately.

\medskip
\noindent{\bf Acknowledgements.} The author is very grateful to A.
Borovik, G. Brightwell and P. Cameron for their helpful comments and
encouragement.

\section{Notation and definitions}

\subsection{Notation}\label{S1.1}

We shall be using the following standard notation.

\begin{enumerate}
\item[$\bullet$]
$\,\N\,$ is the set of positive integers.
\item[$\bullet$]
$\{a_1,a_2,\ldots\}$ is the set with elements $\,a_1,a_2,\ldots\,$.
\item[$\bullet$]
$\#\,A$ denotes the number of elements of a set $\,A\,$.
\item[$\bullet$]
$\supp w\,$ denotes the support of a function $\,w\,$.
\end{enumerate}
If $\,A\,$ is a subset of a real linear space $\,\WC\,$ then
\begin{enumerate}
\item[$\bullet$]
$\ex A$ is the set of extreme points of the set $\,A\,$,
\item[$\bullet$]
$\conv A$ is the convex hull of $\,A\,$.
\end{enumerate}
If $\,\WC\,$ is equipped with a topology $\,\Tf\,$ then
\begin{enumerate}
\item[$\bullet$]
$\cconv A$ is the closure of $\,\conv A\,$.
\end{enumerate}

Let $\,\Ga=\{\Om_1,\Om_2,\ldots\}\,$ be a family of countable sets
$\,\Om_k\,$ which may have nonempty intersection. Denote
$\,\Om:=\cup_k\Om_k\,$, and let
\begin{enumerate}
\item[$\bullet$]
$\SC(\Ga)\,$ be the convex set of nonnegative functions $\,w\,$ on
$\,\Om\,$ such that $\,\sum_{g\in\Om_k}w(g)=1\,$ for all
$\,\Om_k\in\Ga\,$;
\item[$\bullet$]
$\SC^0(\Ga)$ be the convex set of nonnegative functions $\,w\,$ on
$\,\Om\,$ such that  $\,\sum_{g\in\Om_k}w(g)\leq1\,$ for all
$\,\Om_k\in\Ga\,$;
\item[$\bullet$]
$\PC(\Ga)$ be the set of functions $\,w\in\SC(\Ga)\,$ taking only
the values 0 and 1;
\item[$\bullet$]
$\PC^0(\Ga)$  be the set of functions $\,w\in\SC^0(\Ga)\,$ taking
only the values 0 and 1.
\end{enumerate}
If $\,w\in\SC(\Ga)\,$ then the restrictions
$\,\left.w\right|_{\Om_k}\,$ are densities of probability measures
on $\,\Om_k\,$ such that
\begin{equation}\label{1.1}
\left.\mu_i\right|_{\Om_i\cap\Om_j}
=\left.\mu_j\right|_{\Om_i\cap\Om_j}\,, \qquad\forall
i,j=1,2,\ldots\,.
\end{equation}
The other way round, every family of probability measures
$\,\mu_k\,$ on $\,\Om_k\,$ satisfying \eqref{1.1} generates a
function $\,w\in\SC(\Ga)\,$.

If $\,\Om_k\subseteq\Om_j\,$ then all functions $\,w\in\SC(\Ga)\,$
are identically equal to zero on $\,\Om_j\setminus\Om_k\,$. In this
case we can remove the set $\,\Om_j\,$ from $\,\Ga\,$ and all
elements $\,g\in\Om_j\setminus\Om_k\,$ from $\Omega$ without
changing the structure of the sets $\,\SC(\Ga)\,$ and
$\,\PC(\Ga)\,$. Therefore, without loss of generality, we shall
always be assuming that
\begin{equation}\label{1.2}
\Om_k\setminus\Om_j\ne\varnothing\,,\qquad\forall j\ne k\,.
\end{equation}

Given $\,g\in\Om\,$ and a subset $\,\tilde\Om\subseteq\Om\,$, let us
define
\begin{enumerate}
\item[$\bullet$]
$\,\Ga(g):=\{\Om_k\in\Ga\,:\,g\in\Om_k\}\,$,
$\,\Ga(\tilde\Om):=\cup_{g\in\tilde\Om}\,\Ga(g)\,$,
\item[$\bullet$]
$\,\varkappa(g):=\#\,\Ga(g)\,$ and
$\,\varkappa(\tilde\Om):=\sup_{g\in\tilde\Om}\Ga(g)\,$.
\end{enumerate}
Further on, we shall be assuming that
\begin{enumerate}
\item[{\bf(a)}]
$\,\#\,\Ga(g)<\infty\,$ all $\,g\in\Om\,$.
\end{enumerate}
We shall also need the following technical assumption:
\begin{enumerate}
\item[{\bf(a$_1$)}]
if $\,g\ne\tilde g\,$ then $\,\Ga(g)\ne\Ga(\tilde g)\,$ for all
$\,g,\tilde g\in\tilde\Om\,$.
\end{enumerate}

In the following two examples $\,\ka(g)=2\,$ for all $\,g\in\Om\,$
and holds {\bf(a$_1$)} for the whole set $\,\Om=\cup_k\Om_k\,$.

\begin{example}\label{E1.1}
Suppose that one can split the family $\,\Ga\,$ into the union of
two subfamilies $\Ga^+=\{\Om_1^+,\Om_2^+,\ldots\}$ and
$\Ga^-=\{\Om_1^-,\Om_2^-,\ldots\}$ in such a way that
$\,\Om_i^+\cap\Om_j^+=\varnothing\,$,
$\,\Om_i^-\cap\Om_j^-=\varnothing\,$,
$\,\Om=\cup_{i,j}(\Om_i^+\cap\Om_j^-)\,$ and
$\,\#\,(\Om_i^+\cap\Om_j^-)=1\,$ for all $\,i,j=1,2,\ldots\,$ Then
$\,\Om_i^+\,$ and $\,\Om_j^-\,$ may be considered as rows and
columns of an $\,m_+\times m_-$-matrix, where
$\,m^\pm:=\#\,\Ga^\pm\,$. Under this identification, $\,\SC(\Ga)\,$
is the set of doubly stochastic matrices and $\,\PC(\Ga)\,$ is the
set of permutation matrices.
\end{example}

\begin{example}\label{E1.2}
Let $\,\Om=\{g_1,g_2,\ldots,g_n\}\,$, $\,\Om_n:=\{g_n,g_1\}\,$ and
$\,\Om_k:=\{g_k,g_{k+1}\}\,$ for $\,k=1,2,\ldots,n-1\,$. If $\,n\,$
is odd then $\,\PC(\Ga)=\varnothing\,$ and $\,\SC(\Ga)\,$ consists
of one function $\,w\,$ identically equal to $\,\frac12\,$.
\end{example}

\subsection{The set $\,\SC(\Ga)$}\label{S1.2}

It may well happen that $\,\SC(\Ga)=\varnothing\,$. Indeed,
\begin{equation}\label{1.3}
\#\,\Ga\ =\ \sum_{g\in\Om}\ka(g)\,w(g)\ \leq\
\ka(\Om)\,\sum_{g\in\Om}w(g)\,,\qquad\forall w\in\SC(\Ga)\,.
\end{equation}
Therefore $\,\SC(\Ga)=\varnothing\,$ for all finite non-square
matrices (in which the number of rows is not equal to the number of
columns). The equality \eqref{1.3} also implies the following more
general

\begin{lemma}\label{L1.3}
Assume that $\,\Om\,$ coincides with the union of $\,m\,$ sets
$\,\Om_k\in\Ga\,$. If $\,\#\,\Ga>m\,\ka(\Om)\,$ then
$\,\SC(\Ga)=\varnothing\,$.
\end{lemma}

\begin{proof}
Under the conditions of lemma, $\,\sum_{g\in\Om}w(g)\leq m\,$ for
all $\,w\in\SC(\Ga)\,$. If $\,\SC(\Ga)\ne\varnothing\,$ then
$\,\#\,\Ga\leq m\,\ka(\Om)\,$ by virtue of \eqref{1.3}.
\end{proof}

In particular, $\,\SC(\Ga)=\varnothing\,$ whenever
$\,\#\,\Ga=\infty\,$, $\,\ka(\Om)<\infty\,$ and $\,\Om\,$ coincides
with the union of a finite collection of the sets $\,\Om_k\in\Ga\,$.
The following lemma shows that this remains true under the less
restrictive condition {\bf(a)}.

\begin{lemma}\label{L1.4}
Let $\,n\in\N\,$ and $\,w\,$ be a nonnegative function on the union
$\,G_n:=\cup_{k\leq n}\,\Om_k\,$ such that
\begin{equation}\label{1.4}
\sum_{g\in\Om_k}w(g)=1\ \text{for all $k\leq n\ $ and}\
\,\sum_{g\in\Om_k\cap G_n}w(g)\leq1\ \text{for all $k>n$.}
\end{equation}
If $\,\#\,\Ga=\infty\,$ and the condition {\bf(a)} is fulfilled then
\begin{equation}\label{1.5}
\sum_{g\in\Om_j\cap G_n}w(g)\ \to\ 0\,,\qquad j\to\infty\,.
\end{equation}
\end{lemma}

\begin{proof}
Since $\,\sum_{g\in G_n}w(g)\leq n<\infty\,$, for each $\,\ve>0\,$
we can find a finite subset $\,G_{n,\ve}\subseteq G_n\,$ such that
$$
0\ \le\ \sum_{g\in G_n}w(g)-\sum_{g\in G_{n,\ve}}w(g\ )\ <\ve\,.
$$
The condition {\bf(a)} implies that $\,\sup\,\{j\,:\,\Om_j\cap
G_{n,\ve}\}<\infty\,$. Therefore the sum in \eqref{1.5} does not
exceed $\,\ve\,$ for all sufficiently large $\,j\in\N\,$.
\end{proof}

\begin{corollary}\label{C1.5}
Let $\,\#\,\Ga=\infty\,$ and the condition {\bf(a)} be fulfilled. If
$\,\Om\,$ coincides with the union of a finite collection of the
sets $\,\Om_k\in\Ga\,$ then $\,\SC(\Ga)=\varnothing\,$.
\end{corollary}

\subsection{The space $\,\WC\,$}\label{S1.3}

Let $\,\WC\,$ be an arbitrary linear space of real-valued functions
on $\,\Om\,$, which includes $\,\PC^0(\Ga)\,$ and satisfies the
following condition:
\begin{enumerate}
\item[{\bf(w)}]
if $\,w\in\WC\,$ and $\,|\tilde w|\le|w|\,$ then $\,\tilde
w\in\WC\,$.
\end{enumerate}
Since $\,\PC^0(\Ga)\subset\WC\,$, the space $\,\WC\,$ contains all
functions with finite supports. Since $\,\Om\,$ is countable,
$\,\WC\,$ can be thought of as a subspace of the linear space of
infinite real sequences (or a subspace of $\,\R^m\,$, if
$\,\#\,\Om=m<\infty\,$). Let
\begin{enumerate}
\item[$\bullet$]
$\,\SC(\Ga,\WC):=\SC(\Ga)\cap\WC\,$ и
$\,\SC^0(\Ga,\WC):=\SC^0(\Ga)\cap\WC\,$.
\end{enumerate}
Obviously, $\,\PC(\Ga)\subseteq\ex\SC(\Ga,\WC)\,$ and
$\,\PC^0(\Ga)\subseteq\ex\SC^0(\Ga,\WC)\,$.

If $\,\WC\,$ contains the linear space
\begin{equation}\label{1.6}
\WC_1\ :=\ \{\,w\,:\,\sup_k\sum_{g\in\Om_k}|w(g)|<\infty\,\},
\end{equation}
then $\,\SC(\Ga,\WC)=\SC(\Ga)\,$ and
$\,\SC^0(\Ga,\WC)=\SC^0(\Ga)\,$. In particular, this is the case
when $\,\#\,\Om<\infty\,$. If $\,\#\,\Om=\infty\,$ and
$\,\WC_1\not\subset\WC\,$ then the condition $\,w\in\SC(\Ga,\WC)\,$
imposes some restrictions on the behavior of the function
$\,w\in\SC(\Ga)\,$ at infinity.

\begin{remark}\label{R1.6}
In \cite{Sa} I also considered the sets
 $\,\SC(\Ga,\Ga_1)\,$ formed by the nonnegative functions $\,w\,$
such that
$$
\sum_{g\in\Om_k}w(g)\leq1\ \text{for all
}\,\Om_k\in\Ga\quad\text{and}\quad \sum_{g\in\Om_k}w(g)=1\text{ for
all }\,\Om_k\in\Ga_1\,,
$$
where $\,\Ga_1$ is a subset of $\,\Ga\,$. Let $\,\Om'_k\,$ be the
set obtained from $\,\Om_k\,$ by adding one new element  $\,g'_k\,$
that does not belong any other set $\,\Om'_j\,$. If
$\,w\in\SC(\Ga,\Ga_1)\,$ and
$$
w'(g')\ :=\ \begin{cases}w(g')& \text{for all }g'\in\cup_k\Om_k\,,\\
1-\sum_{g\in\Om_k}w(g)&\text{for }g'=g'_k\,,\end{cases}
$$
then $\,w'\in\SC(\Ga')\,$ where $\,\Ga'=\{\Om'_1,\Om'_2,\ldots\}\,$.
Using this observation, one can extend results of this paper to the
sets $\,\SC(\Ga,\Ga_1)\,$.
\end{remark}

\section{Extreme points}

\subsection{The associated graph $G$}\label{S2.1}
Let us join every two elements of $\,\Om=\cup_k\Om_k\,$ lying in the
same set $\,\Om_k\,$ by an edge and consider the graph $\,G\,$
obtained by means of this procedure. In this section we shall
identify subsets of $\,\Om\,$ with subgraphs of $\,G\,$, assuming
that their vertices are adjacent if and only if they are adjacent in
$\,G\,$ (in other words, the subgraphs include all edges joining
their vertices). Then the sets $\,\Om_k\,$ become complete subgraphs
of the graph $\,G\,$. Note that $\,\Om_k\,$ is not necessarily a
maximal complete subgraph (a clique) and that there may be cliques
in $\,G\,$ which do not contain any of the sets $\,\Om_k\,$.

Recall that a path $\,\ga\,$ in $\,G\,$ is s sequence
$\,(g_1,g_2,\dots,g_m)\,$ of vertices $\,g_j\,$ such that
$\,(g_j,g_{j+1})\,$ are distinct graph edges. A path $\,\ga\,$ is
said to be {\it simple} if each its vertex is adjacent to at most
two other vertices of $\,\ga\,$. A path is called a cycle if
$\,g_m=g_1\,$. One says that a cycle is odd (or even) if it contains
an odd (or even) number of vertices.

\begin{definition}\label{D2.1}
We shall call a simple path $\,\ga\in G\,$ {\it primitive} if none
of the sets $\,\Om_k\,$ contain more than two its vertices.
\end{definition}

Note that the graph $\,G\,$ in Example \ref{E1.2} consists of one
primitive cycle with $\,n\,$ vertices.

If a simple path $\,\ga\,$ is not primitive then
\begin{enumerate}
\item[$\bullet$]
either $\,\ga\,$ has exactly three vertices lying in a subgraph
$\,\Om_k\,$,
\item[$\bullet$]
or $\,\ga\,$ contains at least two non-consecutive vertices which
are adjacent in $\,G\,$.
\end{enumerate}
If $\,\ga=(g_1,g_2,\ldots,g_m)$ is a simple path in $\,G\,$ joining
the vertices $\,g_1\,$ and $\,g_m\ne g_1\,$, then the shortest path
of the form $\,(g_1,g_{j_1},g_{j_2}\ldots,g_{j_l},g_m)$ is
primitive. Therefore any two vertices lying in a connected subgraph
$\,G'\,$ can be joined by a primitive path $\,\ga\in G'\,$. In
particular, any path joining two given vertices with the minimal
possible number of vertices is primitive.

\begin{lemma}\label{L2.2}
Let $\,\ga_0=(g_1,g_2,\dots,g_m,g_1)\,$ be a simple cycle. Then
there is a finite collection of cycles
$\,\{\ga'_1,\ga'_2,\ldots\ga'_n\}\,$ such that
\begin{enumerate}
\item[(1)]
each cycle $\,\ga'_i\,$ either is primitive or has exactly three
vertices lying in a subgraph $\,\Om_k\,$;
\item[(2)]
the set of vertices of the cycles $\,\ga'_i\,$ coincides with
$\,\{g_1,g_2,\dots,g_m\}\,$;
\item[(3)]
$\,\ga'_i\,$ and $\,\ga'_{i+1}\,$ have exactly two common vertices
and one common edge;
\item[(4)]
if $\,j\geq2\,$ then $\,\ga'_i\,$ and $\,\ga'_{i+j}\,$ have at most
one common vertex.
\end{enumerate}
\end{lemma}

\begin{proof}
Consider the shortest cycle $\,\ga'_1\,$ of the form
$$
\ga'_1=(g_1,g_2,\dots\,g_j,g_l,g_{l+1},\ldots,g_m,g_1)
$$
which includes the edge $\,(g_1,g_2)\,$. This cycle either is
primitive or has three vertices lying in a subgraph $\,\Om_k\,$.
Denote by $\,\ga_1\,$ the cycle $\,(g_j,g_{j+1},...,g_l,g_j\,)$
obtained from $\,\ga_0\,$ by removing $\,\ga'_1\,$, and enumerate
its vertices such a way that $\,g_l=g_1\,$ and $\,g_j=g_2\,$.
Consecutively applying this procedure to
$\,\ga_0,\ga_1,\ga_2,\ldots\,$, after finitely many steps we obtain
a collection of cycles $\,\ga'_1,\ga'_2,\ldots\ga'_n\,$ satisfying
the conditions (1)--(4).
\end{proof}

\begin{lemma}\label{L2.3}
Let $\,\tilde G\,$ be an arbitrary subgraph of $\,G\,$ satisfying
the condition {\bf(a$_1$)}. If $\,\ka(\tilde G)\leq2\,$ and
$\,\tilde G\,$ does not contain primitive cycles then every simple
cycle $\,\ga_0\in\tilde G\,$ is contained in a subgraph $\,\Om_k\,$.
\end{lemma}

\begin{proof}
Let $\,\ga'_1,\ga'_2,\ldots\ga'_n\,$ be the cycles constructed in
Lemma \ref{L2.2}. If $\,\tilde G\,$ does not contain primitive
cycles then each cycle $\,\ga'_i\,$ lies in a subgraph
$\,\Om_{k_i}\,$. Now the condition (3) of Lemma \ref{L2.2},
{\bf(a$_1$)} and the estimate $\,\ka(\tilde G)\leq2\,$ imply that
$\,\Om_{k_1}=\Om_{k_2}=\dots=\Om_{k_n}\,$.
\end{proof}

\begin{corollary}\label{C2.4}
If a connected subgraph $\,G'\,$ satisfies the conditions of Lemma
{\rm\ref{L2.3}} then every two its vertices are joined by a unique
primitive path.
\end{corollary}

\begin{proof}
Suppose that there are two distinct primitive paths in $\,G'\,$
joining the same two vertices. Then $\,G'\,$ contains a simple cycle
$\,\ga\,$ formed by edges and vertices of these paths. Since no
three consecutive vertices of a primitive path belong to the same
subgraph $\,\Om_k\,$, the cycle $\,\ga\,$ is not contained in any of
these subgraphs. This contradicts to Lemma \ref{L2.3}.
\end{proof}

In the general case, the absence of primitive cycles does not imply
that every two vertices of a connected subgraph $\,G'\,$ are joined
by a unique primitive path.

\begin{example}\label{E2.5}
Let $\,\Om=\{g_0,g_1,\ldots,g_{m+1}\}\,$ and
$\,\Om_k=\{g_0,g_k,g_{k+1}\}\,$, where $\,k=1,\ldots,m\geq2\,$. Then
$\,G\,$ does not contain any primitive cycles but $\,g_1\,$ and
$\,g_{m+1}\,$ are joined by the two distinct primitive paths
$\,(g_1,g_2,\dots,g_m,g_{m+1})\,$ and $\,(g_1,g_0,g_{m+1})\,$. In
this example $\,\ka(g_0)=m\,$ and $\,\ka(g)\leq2\,$ for all other
vertices $\,g\,$. If $\,m>2\,$ then $\,\ka(\Om)>2\,$; if $\,m=2\,$
then {\bf(a$_1$)} is not true because $\,\Ga(g_0)=\Ga(g_2)=\Ga\,$.
\end{example}

\begin{remark}\label{R2.6}
Instead of $\,G\,$, one can consider the so-called intersection
graph where vertices are identified with the sets $\,\Om_k\,$ and
two vertices are adjacent if and only if the intersection of the
corresponding sets is not empty. However, this construction seems to
be less suitable for our purposes because the intersection graph may
contain several distinct edges associated with the same element
$\,g\in\Om\,$.
\end{remark}

\subsection{Necessary conditions}\label{S2.2}

Given $\,w\in\WC\,$, denote $\,G_w:=\supp w\,$ and
$$
\hat w(\tilde G)\ :=\ \inf_{g\in\tilde
G}\min\{w(g),1-w(g)\}\,,\qquad\forall\tilde G\subseteq G_w\,.
$$
If $\,w\in\SC(\Ga)\,$ then
\begin{enumerate}
\item[$\bullet$]
$\,\hat w(\tilde G)\geq0\,$ for all $\,\tilde G\subseteq G_w\,$ and
\item[$\bullet$]
$\,\hat w(G')>0\,$ for every finite connected subgraph $\,G'\,$
containing more than one vertex.
\end{enumerate}

The following two lemmas give necessary conditions on $\,G_w\,$,
which are fulfilled for every extreme point $\,w\,$ of the set
$\,\SC(\Ga,\WC)\,$.

\begin{lemma}\label{L2.7}
If $\,w\in\ex\SC(\Ga,\WC)\,$ then the graph $\,G_w\,$ does not
contain any finite nonempty subgraphs $\,\tilde G\,$ satisfying the
following two conditions:
\begin{enumerate}
\item[(1)]
the intersection of $\,\tilde G\,$ with each subgraph $\,\Om_k\,$
either is empty or contains exactly two distinct vertices;
\item[(2)]
$\,\tilde G\,$ does not contain odd primitive cycles.
\end{enumerate}
\end{lemma}

\begin{proof}
Let $\,w\in\SC(\Ga,\WC)\,$ and $\,\tilde G\,$ be a finite subgraph
of $\,G_w\,$ satisfying the conditions (1) and (2). Let us consider
an arbitrary connected component $\,G'\,$ of the graph $\,\tilde
G\,$. Since the subgraphs $\,\Om_k\,$ are complete, the intersection
$\,G'\cap\Om_k\,$ either is empty or coincides with $\,\tilde
G\cap\Om_k\,$ for each $\,\Om_k\in\Ga\,$. Therefore the subgraph
$\,G'\,$ also satisfies the conditions (1) and (2). In particular,
from (1) it follows that $\,w(g)\in(0,1)\,$ for all $\,g\in G'\,$.
Since $\,G'\,$ is finite, this implies that $\,\hat w(G')>0\,$.

In view of (1), every simple path in $\,G'\,$ is primitive. From
here and the condition (2) it follows that the chromatic number of
the graph $\,G'\,$ equals two. In other words, the vertices of
$\,G'\,$ can be divided into two groups in such a way that no two
vertices from one group are adjacent to each other. Let us denote
$\,\ve:=\hat w(G')\,$ and define functions $\,w_\ve^\pm\,$ as
follows:
\begin{enumerate}
\item[$\bullet$]
$\,w_\ve^\pm(g):=w(g)\,$ for all $\,g\not\in G'\,$;
\item[$\bullet$]
$\,w_\ve^\pm(g):=w(g)\pm\ve\,$ if $g$ belongs to the first group of
vertices;
\item[$\bullet$]
$\,w_\ve^\pm(g):=w(g)\mp\ve\,$ if $g$ belongs to the second group of
vertices.
\end{enumerate}
It follows straight from the definition that
$\,w_\ve^\pm(g)\in\SC(\Ga)\,$. Also, $\,w_\ve^\pm\in\WC\,$ because
$\,w\in\WC\,$ and $\,\#\,\supp(w-w_\ve^\pm)<\infty\,$. Thus
$\,w_\ve^\pm\in\SC(\Ga,\WC)\,$ and
$\,w=\frac12\,(w_\ve^++w_\ve^-)\not\in\ex\SC(\Ga,\WC)\,$.
\end{proof}

\begin{corollary}\label{C2.8}
If $\,w\in\ex\SC(\Ga,\WC)\,$ then each connected component of the
graph $\,G_w\,$ satisfies the condition {\bf(a$_1$)}.
\end{corollary}

\begin{proof}
If $\,g_1\,$ and $\,g_2\,$ lie in the same connected component of
$\,G_w\,$ and $\,\Ga(g_1)=\Ga(g_2)\,$ then the subgraph $\,\tilde
G:=\{g_1,g_2\}\,$ satisfies the conditions (1) and (2) of Lemma
\ref{L2.7}.
\end{proof}

\begin{lemma}\label{L2.9}
If $\,\ka(\Om)\le2\,$ and $\,w\in\ex\SC(\Ga,\WC)\,$ then each
connected component of the graph $\,G_w\,$ either consists of one
vertex or contains a primitive cycle.
\end{lemma}

\begin{proof}
Let $\,w\in\SC(\Ga,\WC)\,$. Suppose that $\,G_w\,$ has a connected
component $\,G'\,$ that includes at least one edge and does not
contain primitive cycles. Then, in view of Corollaries \ref{C2.4}
and \ref{C2.8}, every two vertices of $\,G'\,$ are joined by a
unique primitive path.

Let us fix $\,g_0\in G'\,$ and denote by $\,\GC_n\,$ the set of
vertices in $\,G'\,$ which are joined with $\,g_0\,$ by primitive
paths with $\,n\,$ edges. Since the subgraphs $\,\Om_k\,$ are
complete, for every $k=1,2,\ldots$ there exists $n\geq 0$ such that
$\,\Om_k\subseteq\GC_n\cup\GC_{n+1}\,$. Moreover, if
$\,\Om_k\subseteq\GC_n\cup\GC_{n+1}\,$ then the intersection
$\,\Om_k\cap\GC_n\,$ contains exactly one vertex, which we shall
denote $\,g_{k,n}\,$. Indeed, if $\,\Om_k\cap\GC_n\,$ contained
another vertex $\,g'_{k,n}\,$ then, joining $\,g_0\,$ with
$\,g_{k,n}\,$ and $\,g'_{k,n}\,$ by primitive paths and adding the
edge $\,(g'_{k,n},g_{k,n})\,$, we would obtain a simple cycle not
lying in any of the sets $\,\Om_k\,$. This would contradict to Lemma
\ref{L2.3}.

Since $\,\#\,G'>1\,$, we have $\,w(g_0)\in(0,1)\,$. Let us denote
\begin{equation}\label{2.1}
\ve_0\ :=\
\min\left\{\frac12\,,\;\frac{1-w(g_0)}{2\,w(g_0)}\right\}\,,
\end{equation}
fix an arbitrary $\,\ve\in(0,\ve_0)\,$ and define
$$
\ve_{k,0}:=\ve\,,\qquad \ve_{k,n+1}\ :=\
\ve_{k,n}\,w(g_{k,n})\left(1-w(g_{k,n})\right)^{-1}\,,\qquad
n=1,2,\ldots
$$
Consider the sequences of functions $\,w_{\ve,n}^+\,$ and
$\,w_{\ve,n}^-\,$ defined as follows:
\begin{enumerate}
\item[$\bullet$]
$w_{\ve,0}^\pm(g_0):=(1\pm\ve)\,w(g_0)\,$ and
$\,w_{\ve,0}^\pm(g):=w(g)$ for all $g\ne g_0\,$;
\item[$\bullet$]
$w_{\ve,n+1}^\pm(g):=w_{\ve,n}^\pm(g)\,$ for all $\,g\in\cup_{j\leq
n}\,\GC_j\,$;
\item[$\bullet$]
$w_{\ve,n+1}^\pm(g):=w(g)\,$ for all $\,g\not\in\cup_{j\leq
n+1}\,\GC_j\,$;
\item[$\bullet$]
if $\,\Om_k\subseteq\GC_n\cup\GC_{n+1}\,$ then
$\,w_{\ve,n}^\pm(g_{k,n}):=(1\pm\ve_{k,n})\,w(g_{k,n})\,$ and
\newline
$\,w_{\ve,n+1}^\pm(g):=(1\mp\ve_{k,n+1})w(g)\,$ for all
$\,g\in\Om_k\cap\GC_{n+1}\,$.
\end{enumerate}
Obviously, $\,w(g)=\frac12(w_{\ve,n}^+(g)+w_{\ve,n}^-(g))\,$ for all
$\,g\in\Om\,$. Since $\,w\in\SC(\Ga)\,$, we have
$$
\frac{w(g_{k,n})}{1-w(g_{k,n})}\ \leq\
\frac{1-w(g)}{w(g)}\,,\qquad\forall g\in\Om_k\cap\GC_{n+1}\,.
$$
Using these inequalities, one can easily prove that
\begin{equation}\label{2.2}
\ve_{k,n}\ \leq\
\min\left\{\frac12\,,\,\frac{1-w(g_{k,n})}{2\,w(g_{k,n})}\right\}\,,
\qquad k=1,2,\ldots
\end{equation}
The estimate \eqref{2.2} and identity
$\,w=\frac12(w_\ve^++w_\ve^-)\,$ imply that
\begin{equation}\label{2.3}
0\ \leq\ w_{\ve,n}^\pm(g)\ \leq\ w(g)\,,\qquad\forall g\in
G'\,,\quad\forall n=1,2,\ldots
\end{equation}
If $\,\Om_k\subseteq\GC_n\cup\GC_{n+1}\,$ then
\begin{equation}\label{2.4}
\sum_{g\in\Om_k}w_{\ve,n}^\pm(g) =(1\pm\ve_{k,n})w(g_{k,n}) \ +\
(1\mp\ve_{k,n+1})(1-w(g_{k,n}))\ =\ 1.
\end{equation}

Let $w_\ve^\pm(g):=\lim_{n\to\infty}w_{\ve,n}^\pm(g)$. The condition
{\bf(w)} and inequalities \eqref{2.3} imply that
$\,w_\ve^\pm\in\WC\,$. By \eqref{2.4}, we have
$\,w_\ve^\pm\in\SC(\Ga)\,$. Therefore $\,w_\ve^\pm\in\SC(\Ga,\WC)\,$
and $\,w=\frac12(w_\ve^++w_\ve^-)\,$ is not an extreme point of the
set $\,\SC(\Ga,\WC)\,$.
\end{proof}

\begin{remark}\label{R2.10}
If $\,w\in\PC(\Ga)\,$ then $\,\#\,(\Om_k\cap G_w)=1\,$ for all
$\,\Om_k\in\Ga\,$, which implies that the subgraph $\,G_w\,$ does
not have any edges. Conversely, every edgeless graph containing one
element of each set $\,\Om_k\,$ is the support of a function from
$\,\PC(\Ga)\,$. If $\,w_1,w_2\in\PC(\Ga)\,$ then the subgraph
$\,\tilde G:=\supp(w_1-w_2)\,$ satisfies the conditions (1) and (2)
of Lemma \ref{L2.7}. The converse is not always true; a subgraph
satisfying (1) and (2) may not coincide with the support of the
difference $\,w_1-w_2\,$ with $\,w_1,w_2\in\PC(\Ga)\,$.

For example, let $\,\Om=\{g_1,g_2,\ldots,g_5\}$,
$\,\Om_1=\{g_2,g_3\}\,$, $\,\Om_2=\{g_1,g_3\}\,$,
$\,\Om_3=\{g_1,g_2\}\,$ and $\,\Om_4=\{g_4,g_5\}\,$. Then the
complete subgraph $\,\Om_4\,$ satisfies the conditions (1) and (2).
However, $\,\SC(\Ga)\,$ consists of functions $\,w\,$ such that
$\,w(g_1)=w(g_2)=w(g_3)=\frac12\,$ and
$\,1-w(g_5)=w(g_4)\in[0,1]\,$. Therefore $\,\PC(\Ga)=\varnothing\,$.
\end{remark}

\subsection{Necessary and sufficient conditions}\label{S2.3}

The following theorem is the main result of this section.

\begin{theorem}\label{T2.11}
Assume that $\,\ka(\Om)\leq2\,$. Then $\,w\in\ex\SC(\Ga,\WC)\,$ if
and only if $\,G_w\cap\Om_k\ne\varnothing\,$ for all
$\,\Om_k\in\Ga\,$ and $\,w\,$ satisfies the following two
conditions:
\begin{enumerate}
\item[(1)]
each connected component of the graph $\,G_w\,$
\begin{enumerate}
\item[(1$_1$)]
either consists of one isolated vertex $\,g\,$,
\item[(1$_2$)]
or coincides with an odd primitive cycle $\,\ga\,$,
\end{enumerate}
\item[(2)]
$\,w(g)=1\,$ in the case {\rm (1$_1$)}, and
$\,\left.w\right|_\ga\equiv\frac12\,$ in the case {\rm (1$_2$)}.
\end{enumerate}
\end{theorem}

\begin{proof} If (1) holds and $\,G_w\cap\Om_k\ne\varnothing\,$
then
\begin{enumerate}
\item[(1$'$)]
for each $\,\Om_k\in\Ga\,$ the intersection $\,\Om_k\cap G_w\,$
consists of
\begin{enumerate}
\item[(1$_1'$)]
either one isolated vertex $\,g\in G_w\,$,
\item[(1$_2'$)]
or two consecutive vertices of a primitive cycle $\,\ga\in G_w\,$.
\end{enumerate}
\end{enumerate}
Indeed, since the graph $\,\Om_k\,$ is complete, the intersection
$\,\Om_k\cap G_w\,$ lies in a connected component $\,G'\,$ of the
graph $\,G_w\,$. If $\,G'\,$ is an isolated vertex then (1$_1'$) is
true. If $\,G'\,$ coincides with a primitive cycle then the
intersection $\,\Om_k\cap G_w\,$ cannot contain more than two
vertices because the cycle is primitive. On the other hand,
$\,\Om_k\cap G_w\,$ cannot consists of one vertex because
$\,\ka(\Om)\leq2\,$. Therefore (1$_2'$) holds.

\smallskip\noindent 1. \
Let $\,G_w\cap\Om_k\ne\varnothing\,$ for all $\,\Om_k\in\Ga\,$ and
the conditions (1) and (2) be fulfilled.

The set of vertices of an odd primitive cycle can be represented as
the union of three disjoint subsets such that
$\,\Ga(g)\cap\Ga(g')=\varnothing\,$ for every pair of elements
$\,g,g'\,$ lying in the same subset. Therefore $\,w\,$ coincides
with $\,\frac12\,(w_1+w_2+w_3)\,$ where
$\,w_i\in\PC^0(\Ga)\subset\WC\,$, $\,i=1,2,3\,$. This observation
and the conditions (1$'$), (2) imply that $\,w\in\SC(\Ga,\WC)\,$.

Assume that $\,w_\pm\in\SC(\Ga,\WC)\,$ and
$\,w=\frac12\,(w_++w_-)\,$. Then $\,\supp w_\pm\subseteq G_w\,$ and
$\,w_\pm(g)=1\,$ at all isolated vertices $\,g\in G_w\,$. If
$\,\ga\,$ is a cycle in $\,G_w\,$ then $\,w_\pm(g)+w_\pm(g')=1\,$
for each pair of consecutive vertices $\,g,g'\in\ga\,$. Since
$\,\ga\,$ is odd, these equalities imply that $\,w_\pm(g)=\frac12\,$
at all vertices $\,g\in\ga\,$. Thus $\,w_\pm=w\,$ and, consequently,
$\,w\in\ex\SC(\Ga,\WC)\,$.

\smallskip\noindent 2. \
Let $\,w\in\SC(\Ga,\WC)\,$ and (1) be fulfilled. Then
$\,G_w\cap\Om_k\ne\varnothing\,$ for all $\,\Om_k\in\Ga\,$. In the
case (1$_1'$)  $\,w(g)=1\,$; in the case (1$_2'$) $\,w(g)=\frac12\,$
for all $\,g\in\ga\,$ because the cycle $\,\ga\,$ is odd. Thus (2)
follows from the inclusion $\,w\in\SC(\Ga,\WC)\,$ and (1).

\smallskip\noindent 3. \
It remains to show that (1) holds for all $\,w\in\ex\SC(\Ga,\WC)\,$.
Assume that $\,w\in\SC(\Ga,\WC)\,$ and consider an arbitrary
connected component $\,G'\,$ of the graph $\,G_w\,$ containing more
than one vertex. We are going to prove that
$\,w\not\in\ex\SC(\Ga,\WC)\,$ unless $\,G'\,$ coincides with an odd
primitive cycle. In view of Corollary \ref{C2.8} and Lemma
\ref{L2.9}, we can assume without loss of generality that $\,G'\,$
satisfies {\bf(a$_1$)} and contains at least one primitive cycle.

\smallskip\noindent 4. \
If there is an even primitive cycle $\,\ga\subset G'\,$ then the
estimate $\,\ka(\Om)\leq2\,$ implies that the subgraph $\,\tilde
G=\ga\,$ satisfies the conditions of Lemma \ref{L2.7} and,
consequently, $\,w\not\in\ex\SC(\Ga,\WC)\,$. Therefore we shall be
assuming that $\,G'\,$ does not contain even primitive cycles.

\smallskip\noindent 5. \
Let $\,\ga=(g_1,g_2,\dots,g_l,g_1)\,$ be an odd primitive cycle in
$\,G'\,$. Since $\,G_w\,$ satisfies {\bf(a$_1$)}, for each pair of
consecutive vertices $\,g_i,g_{i+1}\in\ga\,$ there exists a unique
set $\,\Om_{k_i}\,$ containing these vertices (we take
$\,g_{l+1}:=g_1\,$). Denote by $\,G_\ga\,$ the graph generated by
the family of sets $\,\Ga(\ga)=\{\Om_{k_1},\ldots,\Om_{k_l}\}\,$.

Let $\,G'\setminus\ga\ne\varnothing\,$. Since the graph $\,G'\,$ is
connected, the estimate $\,\ka(\Om)\leq2\,$ implies that at least
one of the intersections $\,\Om_{k_i}\cap G'\,$ contains a vertex
that does not belong to $\,\ga\,$. Assume, for the sake of
definiteness, that $\,\Om_{k_1}\cap G'\,$ contains a vertex
$\,g_0\not\in\ga\,$.

\smallskip\noindent 5(a). \
If $\,g_0\,$ belongs to a set $\,\Om_{k_j}\,$ with $\,j>1\,$ then
one of the primitive cycles $\,(g_0,g_2,\dots,g_j,g_0)\,$ and
$\,(g_0,g_{j+1},\dots,g_l,g_1,g_0)\,$ is even. If there exists a
primitive path $\,(g_0,\tilde g_1,\dots,\tilde g_m)\,$ in $\,G'\,$
such that $\,\tilde g_m\in\Om_{k_j}\,$ and $\,\tilde g_i\not\in
G_\ga\,$ for all $\,i=1,\ldots,m-1\,$ then
\begin{enumerate}
\item[$\bullet$]
in the case $\,j>1\,$ either $\,(g_0,\tilde g_1,\dots,\tilde
g_m,g_{j+1},g_{j+2},\dots,g_l,g_1,g_0)\,$ or  $\;(g_0,\tilde
g_1,\dots,\tilde g_m,g_j, g_{j-1},\dots,g_2,g_0)\,$ is an even
primitive cycle;
\item[$\bullet$]
in the case $\,j=1\,$ one of the primitive cycles $\,(g_0,\tilde
g_1,\dots,\tilde g_m,g_0)\,$ and $\,(g_0,\tilde g_1,\dots,\tilde
g_m,g_2,\dots,g_l,g_1,g_0)\,$ is even.
\end{enumerate}
Therefore the absence of even primitive cycles implies the following
condition:
\begin{enumerate}
\item[(5a)]
if the first edge of the path $\,(g_0,\tilde g_1,\ldots,\tilde
g_m)\,$ is not contained in the subgraph $\,\Om_{k_1}\,$ then none
of its vertices $\,\tilde g_i\,$ belongs to $\,G_\ga\,$.
\end{enumerate}

Denote by $\,G''\,$ the subgraph of $\,G'\,$ formed by the vertex
$\,g_0\,$ and all the vertices $\,g\,$ which are joined with
$\,g_0\,$ by such paths. By (5a), we have $\,G_\ga\cap
G''=\{g_0\}\,$. If $\,\Om_k\not\in\Ga(\ga)\,$ and $\,\Om_k\cap
G''\ne\varnothing\,$ then $\,\Om_k\subseteq G''\,$ because the
subgraph $\,\Om_k\,$ is complete. Therefore either $\,G''=\{g_0\}\,$
or $\,G''\,$ coincides with the graph generated by the family of
sets $\,\Ga''=\{\Om_{n_1},\Om_{n_2},\ldots\}\subset\Ga\,$ such that
$\,\Om_{n_i}\not\in\Ga(\ga)\,$ and $\,\Om_{n_i}\cap
G''\ne\varnothing\,$.

\smallskip\noindent 5(b). \
Assume that $\,G''=\{g_0\}\,$. Then $\,g_0\,$ belongs only to the
set $\,\Om_{k_1}\,$. Let
\begin{equation}\label{2.5}
 \ve_1\ :=\ \hat w(\{g_0,g_1,\ldots,g_l\})
\end{equation}
and $\,\ve\in(0,\ve_1)\,$. Consider the functions $\,w_\ve^+\,$ and
$\,w_\ve^-\,$ defined as follows:
\begin{enumerate}
\item[$\bullet$]
$\,w_\ve^\pm(g)=w(g)\,$ for all
$\,g\not\in\{g_0,g_1,\ldots,g_l\}\,$;
\item[$\bullet$]
$\,w_\ve^\pm(g_0)=w(g_0)\pm\ve\,$;
\item[$\bullet$]
$\,w_\ve^\pm(g_1)=w(g_1)\mp\ve/2\,$;
\item[$\bullet$]
$\,w_\ve^\pm(g_i)=w(g_1)\pm(-1)^{i-1}\,\ve/2\,$ for all
$\,i=2,\ldots,l\,$.
\end{enumerate}
Since $\,\ka(\Om)\leq2\,$, the vertex $\,g_1\,$ belongs only to the
sets $\,\Om_{k_1}\,$ and $\,\Om_{k_l}\,$, and each vertex $\,g_i\,$
with $\,i\geq 2\,$ belongs only to the sets $\,\Om_{k_{i-1}}\,$ and
$\,\Om_{k_i}\,$. Therefore $\,w_\ve^\pm\in\SC(\Ga,\WC)\,$ and
$\,w=\frac12\,(w_\ve^++w_\ve^-)\not\in\ex\SC(\Ga,\WC)\,$.

\smallskip\noindent 5(c). \
Assume that $\,G''\,$ has more than one vertex and does not contain
primitive cycles. Let $\,0<\ve<\min\{\ve_0,\ve_1\}\,$ where
$\,\ve_0\,$ and $\,\ve_1\,$ are defined in \eqref{2.1} and
\eqref{2.5}. Consider the functions $\,w_\ve^+\,$ and $\,w_\ve^-\,$,
defined as follows:
\begin{enumerate}
\item[$\bullet$]
if $\,g\not\in\ga\,$ and $\,g\not\in G''\,$ then
$\,w_\ve^\pm(g):=w(g)\,$;
\item[$\bullet$]
if $\,g\in\ga\,$ then $\,w_\ve^\pm(g)\,$ is defined as in the part
5(b);
\item[$\bullet$]
if $\,g\in G''\,$ then $\,w_\ve^\pm(g)\,$ is defined as in the proof
of Lemma \ref{L2.9}.
\end{enumerate}
Since $\,G_\ga\cap G''=\{g_0\}\,$, the inequalities
$\,\ka(\Om)\leq2\,$ and \eqref{2.3} imply that
$\,w_\ve^\pm\in\SC(\Ga,\WC)\,$. Thus
$\,w=\frac12\,(w_\ve^++w_\ve^-)\not\in\ex\SC(\Ga,\WC)\,$.

\smallskip\noindent 5(d). \
Finally, let us assume that $\,G''\,$ contains an odd primitive
cycle $\,\ga'\,$. In view of (5a), $\,g_0\,$ does not belong to
$\,\ga'\,$. Let us consider the shortest primitive path
$\,\ga''=(g_0,g''_1,\dots,g''_m,g')\,$ joining $\,g_0\,$ with
$\,\ga'\,$, and enumerate vertices of the cycle
$\,\ga'=(g'_1,g'_2,\dots,g'_n,g'_1)\,$ in such a way that
$\,g'_1=g'\,$. Since $\,\ka(\Om)\leq2\,$, either the three vertices
$\,g''_m,g'_n,g'_1\,$ or the three vertices $\,g''_m,g'_1,g'_2\,$
belong to the same set $\,\Om_k\,$. Assume, for the sake of
definiteness, that the latter is true.

Let
$$
\ve_2\ :=\ \hat
w\{g_0,g_1,\ldots,g_l,g'_1,\ldots,g'_n,g''_1,\ldots,g''_m\}\,,\qquad
\ve\in(0,\ve_2)\,,
$$
and $\,w_\ve^\pm\,$ be the functions defined as follows:
\begin{enumerate}
\item[$\bullet$]
if $\,g\not\in\ga\cup\ga'\cup\ga''\,$ then $\,w_\ve^\pm(g):=w(g)\,$;
\item[$\bullet$]
if $\,g\in\ga\,$ or $\,g=g_0\,$ then $\,w_\ve^\pm(g)\,$ is defined
as in the part 6(b);
\item[$\bullet$]
$\,w_\ve^\pm(g''_i)=w(g''_i)\pm(-1)^i\,\ve\,$ for all
$\,i=1,\ldots,m\,$;
\item[$\bullet$]
$\,w_\ve^\pm(g'_1)=w(g'_1)\pm(-1)^{m+1}\,\ve/2\,$;
\item[$\bullet$]
$\,w_\ve^\pm(g'_i)=w(g'_i)\pm(-1)^{m+i-1}\,\ve/2\,$ for all
$\,i=2,\ldots,n\,$.
\end{enumerate}
Since $\,\ka(\Om)\leq2\,$, we have $\,w_\ve^\pm\in\SC(\Ga,\WC)\,$.
Therefore $\,w=\frac12\,(w_\ve^++w_\ve^-)\,$ is not an extreme point
of $\,\SC(\Ga,\WC)\,$.
\end{proof}

Theorem \ref{T2.11} immediately implies

\begin{corollary}\label{C2.12}
Let $\,\ka(\Om)\leq2\,$. Then
\begin{enumerate}
\item[(1)]
$\,w(g)\in\{0,\frac12,1\}\,$ for all $\,w\in\ex\SC(\Ga,\WC)\,$ and
$\,g\in\Om\,$;
\item[(2)]
$\,\ex\SC(\Ga,\WC)=\PC(\Ga)\,$ whenever $\,G\,$ does not contain odd
primitive cycles.
\end{enumerate}
\end{corollary}

\begin{remark}\label{R2.13}
If $\,\ka(\Om)\leq2\,$ and $\,G\,$ does not contain odd primitive
cycles then one can split $\,\Ga\,$ into the union of two disjoint
subsets $\,\Ga^+=\{\Om_1^+,\Om_2^+,\ldots\}\,$ and
$\,\Ga^-=\{\Om_1^-,\Om_2^-,\ldots\}\,$ in such a way that
$\,\Om_i^+\cap\Om_j^+=\varnothing\,$ and
$\,\Om_i^-\cap\Om_j^-=\varnothing\,$ for all $\,i,j=1,2,\ldots\,$
($\,\Om_i\,$ and $\,\Om_j\,$ are included into the same set
$\,\Ga^\pm\,$ if every primitive path joining $\,\Om_i\,$ and
$\,\Om_j\,$ has an odd number of edges). If, in addition,
$\,\Om=\cup_{i,j}(\Om_i^+\cap\Om_j^-)\,$ and
$\,\#\,(\Om_i^+\cap\Om_j^-)=1\,$ for all $\,i,j=1,2,\ldots\,$ then
$\,\SC(\Ga)\,$ can be thought of as a set of doubly stochastic
matrices (see Example \ref{E1.1}).
\end{remark}

\section{Closed convex hull of extreme points}

\subsection{Topologies on $\,\WC'\,$}\label{S3.1}

Let $\,\WC'\,$ be the space of real-valued functions $\,w'\,$ on
$\,\Om\,$ such that $\,\sum_{g\in\Om}|w(g)\,w'(g)|<\infty\,$ for all
$\,w\in\WC\,$. Further on we shall identify functions
$\,w'\in\WC'\,$ with the corresponding linear functionals
$$
w\ \to\ \langle w,w'\rangle\ :=\ \sum_{g\in\Om}w(g)\,w'(g)
$$
on the space $\,\WC\,$.

Let $\,\Tf\,$ be a locally convex topology on $\,\WC\,$ satisfying
the following conditions:
\begin{enumerate}
\item[{\bf(w$_1$)}]
the topological dual $\,\WC^*\,$ is a subspace of $\,\WC'\,$;
\item[{\bf(w$_2$)}]
$\,\WC^*\,$ contains all the functionals
$\,w\to\sum_{g\in\Om_k}w(g)\,$, $\,k=1,2,\ldots$
\end{enumerate}
From the condition {\bf(w$_2$)} it follows that
$\,\cconv\ex\SC(\Ga,\WC)\subseteq\SC(\Ga,\WC)\,$ and
$\,\cconv\ex\SC^0(\Ga,\WC)\subseteq\SC^0(\Ga,\WC)\,$.

Denote by $\,\Tf_0\,$ the topology of pointwise convergence on
$\,\WC\,$. If {\bf(w$_1$)} holds and $\,\WC^*\,$ consists of
functions with finite supports then $\,\Tf=\Tf_0\,$. The Tikhonov
theorem and Fatou lemma imply that the set $\,\SC^0(\Ga)\,$ is
$\,\Tf_0$-compact. Therefore, by the Krein--Milman theorem,
$\,\SC^0(\Ga)\,$ coincides with the $\,\Tf_0$-closure of the set
$\,\conv\ex\SC^0(\Ga)\,$.

If $\,\#\,\Om<\infty\,$ then $\,\dim\WC<\infty\,$, $\,\Tf=\Tf_0\,$,
the set $\,\SC(\Ga,\WC)\,$ is a compact convex polytope and,
consequently, $\,\SC(\Ga,\WC)=\conv\ex\SC(\Ga,\WC)\,$. If at least
one of the sets $\,\Om_k\,$ is infinite then it may well happen that
the set $\,\SC(\Ga,\WC)\,$ is not compact in any locally convex
topology on $\,\WC\,$. Indeed, if the linear functional
$\,w\to\sum_{g\in\Om_k}w(g)\,$ is not continuous then, as a rule,
the set $\,\SC(\Ga,\WC)\,$ is not closed (there are exceptions, for
instance when $\,\SC(\Ga,\WC)=\varnothing\,$, but such exotic
examples hardly deserve serious consideration). On the other hand,
if the functional $\,w\to\sum_{g\in\Om_k}w(g)\,$ is continuous and
$\,G_k=\{g_1,g_2,\ldots\}\,$ then $\,\SC(\Ga,\WC)\,$ may be compact
only under the very restrictive assumption that
$$
\sup\limits_{w\in\SC(\Ga,\WC)}\sum_{i>j}
|w(g_i)|\underset{j\to\infty}\to0
$$
(this follows from Theorem 1.2 in \cite{Sa} which is proved in the
same way as Theorem 2.4 in the second chapter of \cite{Ru}).

In all known to us examples either
$\,\cconv\ex\SC(\Ga,\WC)=\SC(\Ga,\WC)\,$ or at least one of the
conditions {\bf(w$_1$)} and {\bf(w$_2$)} is not satisfied (see, for
instance, \cite{Is} or Remark \ref{R3.7} in the end of the section).
It is quite possible that these conditions are sufficient. However,
we can prove that $\,\cconv\ex\SC(\Ga,\WC)=\SC(\Ga,\WC)\,$ only
under the following additional assumption:
\begin{enumerate}
\item[{\bf(a$_2$)}]
there exists $\,m\in\N\,$ such that
\begin{enumerate}
\item[{\bf(a$_{21}$)}]
either $\,\Om=\cup_{j=1}^m\Om_j\,$,
\item[{\bf(a$_{22}$)}]
or none of the sets $\,\Om_k\in\{\Om_{m+1},\Om_{m+2},\ldots\}\,$
lies in the union of a finite collection of other sets
$\,\Om_j\in\Ga\,$.
\end{enumerate}
\end{enumerate}
In particular, the condition {\bf(a$_2$)} is satisfies if
$\,\#\,\Ga<\infty\,$, or if the number of finite sets $\,\Om_k\,$ is
finite and $\,\#\,(\Om_i\cap\Om_j)<\infty\,$ for all $\,i\ne j\,$.

\subsection{An extension lemma}\label{S3.2}

Let $\,G_n:=\cup_{k\leq n}\,\Om_k\,$, and let
\begin{enumerate}
\item[$\bullet$]
$\,T_0$ be the operator of extension by zero from $\,G_n\,$ to
$\,\Om\,$,
\item[$\bullet$]
$\,\SC_n^0(\Ga,\WC)\,$ be the convex set of nonnegative functions
$\,w\,$ on $\,G_n\,$ satisfying \eqref{1.4} and such that
$\,T_0w\in\WC\,$.
\end{enumerate}

The role of the condition {\bf(a$_{22}$)} is clarified in the
following lemma.

\begin{lemma}\label{L3.1}
If $\,\#\,\Ga=\infty\,$ and the conditions {\bf(w)}, {\bf(a)},
{\bf(a$_{22}$)} are fulfilled then for each $\,n>m\,$ there exists a
(nonlinear) extension operator $\,T_n\,$ from $\,G_n\,$ to $\,\Om\,$
such that
\begin{enumerate}
\item[(1)]
$\,T_n:\SC_n^0(\Ga,\WC)\to\SC(\Ga,\WC)\,$,
\item[(2)]
$\,T_n:\ex\SC_n^0(\Ga,\WC)\to\ex\SC(\Ga,\WC)\,$,
\item[(3)]
$\,\sup\limits_{w\in\SC_n^0(\Ga,\WC)} \langle
T_0w-T_nw,w'\rangle\underset{n\to\infty}\to0\,$ for all
$\,w'\in\WC'\,$.
\end{enumerate}
\end{lemma}

\begin{proof}\

\smallskip\noindent 1. \
Consider an arbitrary function $\,w\in\SC_n^0(\Ga,\WC)\,$. Let
$\,\Ga_1$ be the family of all sets $\,\Om_k\in\Ga\,$ such that
$\,\sum_{g\in\Om_k\cap G_n}w(g)=1\,$, $\,\GC_1$ be the union of
these sets, and $\,w_1$ be the extension of $\,w\,$ by zero to
$\,\GC_1\,$. By Lemma\ref{L1.4},
$$
\de_j(w_1)\ :=\ \sum_{g\in\Om_j\cap\GC_1}w_1(g)\
\underset{j\to\infty}\to\ 0.
$$
In particular, this implies that $\,\#\,\Ga_1<\infty\,$.

Let $\,k_1:=\min\{k\,\,:\Om_k\not\in\Ga_1\}\,$, $\,\Ga_2\,$ be the
family of sets obtained from $\,\Ga_1\,$ by adding the set
$\,\Om_{k_1}\,$, and $\,\GC_2$ be the union of sets
$\,\Om_k\in\Ga_2\,$. In view of the condition {\bf(a$_{22}$)}, we
have $\,\Om_{k_1}\setminus G_j\ne\varnothing\,$ for all $\,j>k_1\,$.
Let us choose an index $\,j\,$ such that
$\,\de_i(w_1)<\de_{k_1}(w_1)\,$ for all $\,i>j\,$ and fix an
arbitrary element $\,g_1\in\Om_{k_1}\setminus G_j\,$. Let $\,w_2\,$
be the function on $\,\GC_2\,$ defined by the equalities
$\,\left.w_2\right|_{\GC_1}:=\left.w_1\right|_{\GC_1}\,$,
$\,w_2(g_1):=1-\de_{k_1}(w_1)\,$ and $\,w_2(g'):=0\,$ at all other
vertices $\,g'\,$. Then $\,\sum_{g\in\Om_k}w_2(g)=1\,$ for all
$\,\Om_k\in\Ga_2\,$, $\,\sum_{g\in\Om_k\cap\GC_2}w_2(g)<1\,$ for all
$\,\Om_k\in\Ga\setminus\Ga_2\,$ and, by Lemma \ref{L1.4},
$$
\de_j(w_2)\ :=\ \sum_{g\in\Om_j\cap\GC_2}w_2(g)\
\underset{j\to\infty}\to\ 0.
$$

Let $\,k_2:=\min\{k\,\,:\Om_k\not\in\Ga_2\}\,$,
$\,\Ga_3:=\Ga_2\cup\{\Om_{k_2}\}\,$ and $\,\GC_3$ be the union of
sets $\,\Om_k\in\Ga_3\,$. Let us define
$\,\Ga'_2:=\Ga_1\cup\Ga(g_1)\,$, choose $\,j\in\N\,$ so large that
$\,\de_i(w_2)<\de_{k_2}(w_2)\,$ for all $\,i>j\,$ and fix an element
$\,g_2\in\Om_{k_2}\setminus G_j\,$ such that
$\,\Ga(g_2)\cap\Ga'_2=\varnothing\,$ (in view of {\bf(a)} and
{\bf(a$_{22}$)} such an element does exist). Let $\,w_3\,$ be the
function on $\,\GC_3\,$ defined by the equalities
$\,\left.w_3\right|_{\GC_2}:=\left.w_2\right|_{\GC_2}\,$,
$\,w_3(g_2):=1-\de_{k_2}(w_2)\,$ and $\,w_3(g'):=0\,$ at all other
vertices $\,g'\,$. Then
$$
\de_j(w_3)\ :=\ \sum_{g\in\Om_j\cap\GC_3}w_3(g)\
\underset{j\to\infty}\to\ 0.
$$

Now, let $\,k_3:=\min\{k\,\,:\Om_k\not\in\Ga_3\}\,$ and $\,\GC_4$ be
the union of all sets $\,\Om_k\in\Ga_4:=\Ga_3\cup\{\Om_{k_3}\}\,$.
Let us define $\,\Ga'_3:=\Ga'_2\cup\Ga(g_2)\,$, choose $\,j\in\N\,$
such that $\,\de_i(w_3)<\de_{k_3}(w_3)\,$ for all $\,i>j\,$, fix an
element $\,g_3\in\Om_{k_3}\setminus G_j\,$ such that
$\,\Ga(g_3)\cap\Ga'_3=\varnothing\,$, and consider the function
$\,w_4\,$ on $\,\GC_4\,$ defined by the equalities
$\,\left.w_4\right|_{\GC_3}=\left.w_3\right|_{\GC_3}\,$,
$\,w_4(g_3)=1-\de_{k_3}(w_3)\,$ and $\,w_4(g')=0\,$ at all other
vertices $\,g'\,$.

Iterating this procedure, we obtain sequences of families
$\,\Ga_j\subset\Ga\,$ of sets $\,\Om_k\,$, their unions
$\,\GC_j\subset\Om\,$, elements $\,g_j\in\GC_{j+1}\setminus\GC_j\,$
and functions $\,w_j\,$ on the sets $\,\GC_j\,$ such that
\begin{enumerate}
\item[(c$_1$)]
$\,\Ga_j\subset\Ga_{j+1}\,$ and $\,\cup_{j=1}^\infty\Ga_j=\Ga\,$;
\item[(c$_2$)]
$\,\left.w_{j+1}\right|_{\GC_j}=\left.w_j\right|_{\GC_j}\,$ and
$\,\supp(w_{j+1}-w_j)=\{g_j\}\subset\Om_{k_j}\,$;
\item[(c$_3$)]
$\,\sum_{g\in\Om_k\cap\GC_j}w_j(g)\leq1\,$ for all $\,\Om_k\in\Ga\,$
and $\,\sum_{g\in\Om_k}w_j(g)=1\,$ for all $\,\Om_k\in\Ga_j\,$;
\item[(c$_4$)]
for each $\,j\in\N\,$ there exists at most one index $\,i<j\,$ such
that $\,\Ga(g_j)\cap\Ga(g_i)\ne\varnothing\,$.
\end{enumerate}

Let $\,T_n\,$ be the operator defined by the equality
$\,T_nw:=\lim_{j\to\infty}T_0w_j\,$ where the limit is taken in the
topology of pointwise convergence $\,\Tf_0\,$. Since $\,0\leq
T_nw\leq1\,$, from (c$_1$)--(c$_3$) it follows that
$\,T_nw\in\SC(\Ga)\,$.

The condition (c$_2$) also implies that
\begin{enumerate}
\item[(c$_2'$)]
$\supp(T_nw-T_0w)=\ga\,$ where $\ga:=\{g_1,g_2,\ldots\}\,$.
\end{enumerate}
Let us define subsets $\,\ga'\,$ and $\,\ga''\,$ of $\,\ga\,$ as
follows:
\begin{enumerate}
\item[$\bullet$]
$g_1\in\ga'\,$ and $\,g_j\in\ga'\,$ whenever
$\,\Ga(g_j)\cap\Ga(g_i)=\varnothing\,$ for all $\,i<j\,$;
\item[$\bullet$]
if $\,\Ga(g_j)\cap\Ga(g_i)\ne\varnothing\,$ for some $\,i<j\,$ then
\begin{enumerate}
\item[$\centerdot$]
$\,g_j\in\ga'\,$ in the case where $\,g_i\in\ga''\,$,
\item[$\centerdot$]
$\,g_j\in\ga''\,$ in the case where $\,g_i\in\ga'\,$.
\end{enumerate}
\end{enumerate}
In view of (c$_4$), $\,\ga=\ga'\cup\ga''\,$,
$\,\ga'\cap\ga''=\varnothing\,$ and
$\,\Ga(g')\cap\Ga(g'')=\varnothing\,$ for every pair of distinct
elements $\,g',g''\in\ga'\,$ and every pair of distinct elements
$\,g',g''\in\ga''\,$. Therefore the characteristic functions
$\,\chi'\,$ and $\,\chi''\,$ of the sets $\,\ga'\,$ and $\,\ga''\,$
belong to $\,\PC^0(\Ga)\,$. Since $\,T_nw-T_0w\leq1\,$, we have
$\,T_nw-T_0w\leq\chi'+\chi''\,$. This estimate and the condition
{\bf(w)} imply that $\,T_nw\in\WC\,$. Thus
$\,T_nw\in\SC(\Ga,\WC)\,$, that is, $\,T_n\,$ satisfies (1).

\smallskip\noindent 2. \
Suppose that $\,w\in\ex\SC_n^0(\Ga,\WC)\,$ and
$\,T_nw\not\in\ex\SC(\Ga,\WC)\,$. Then there exists a function
$\,\tilde w\not\equiv0\,$ on the set $\,\Om\,$ such that
$\,T_nw\pm\tilde w\in\SC(\Ga,\WC)\,$. The condition (c$_2'$) and the
inclusions $\,T_nw\in\SC(\Ga)\,$ and $\,T_nw\pm\tilde
w\in\SC(\Ga)\,$ imply that
$$
\supp\tilde w\ \subseteq\ \supp T_nw\ \subseteq\
G_n\cup\{g_1,g_2,\ldots\}
$$
and $\,\sum_{g\in\Om_j\cap G_n}|\tilde w(g)|\leq1\,$ for all
$\,\Om_k\in\Ga\,$.

Denote by $\,w_\ve^\pm\,$ the restrictions of the nonnegative
functions $\,T_nw\pm\ve\tilde w\,$ to the set $\,G_n\,$. By Lemma
\ref{L1.4},
$$
\de\ :=\ \sup_{\Om_j\not\in\Ga_1}\sum_{g\in\Om_j\cap\GC_1}w(g)\ <\
1\,.
$$
Therefore $\,w_\ve^\pm\in\SC_n^0(\Ga,\WC)\,$ for all
$\,\ve<1-\de\,$. Since $\,w_\ve^++w_\ve^-=w\,$, from here and the
inclusion $\,w\in\ex\SC_n^0(\Ga,\WC)\,$ it follows that
$\,\left.\tilde w\right|_{\GC_1}\equiv0\,$. Now, using (c$_2$) and
the equalities
$$
w_{j+1}(g_j)\ =\ 1-\de_{k_j}(w_j)\ =\
1\,-\sum_{g\in\Om_{k_j}\setminus\{g_j\}}w_j(g)\,,
$$
by induction in $\,j\,$ we obtain $\,\tilde w(g_j)=0\,$ for all
$\,g_j\,$. Thus $\,\tilde w\equiv0\,$. This contradiction proves
(2).

\smallskip\noindent 3. \
Let $\,w'\in\WC'\,$. Suppose that
$\,\sup\limits_{w\in\SC_n^0(\Ga,\WC)}\langle T_0w-T_nw,w'\rangle\,$
does not converge to zero. Then there exist $\,\de>0\,$ and a
sequence of functions $\,w_n\in\SC_n^0(\Ga,\WC)\,$ such that
$\,\langle T_0w_n-T_nw_n,w'\rangle\geq\de\,$. If $\,\tilde
w_n:=T_nw_n-T_0w_n\,$ then $\,\tilde w_n\in\SC^0(\Ga,\WC)\,$,
$\,\supp\tilde w_n=\{g_1^n,g_2^n,\ldots\}\subset\Om\setminus G_n\,$
and
$$
\sum_{j=1}^\infty|\tilde w_n(g_j^n)\,w'(g_j^n)|\ \geq\ \langle
T_nw_n-T_0w_n,w'\rangle\ \geq\ \de\,,
$$
where $\,\{g_1^n,g_2^n,\ldots\}\,$ are the sets of vertices
associated with functions $\,w_n\,$ (see the first part of the
proof).

Let us consider arbitrary finite subsets
$\,H_n\subset\{g_1^n,g_2^n,\ldots\}\,$ such that
$$
\sum_{g\in H_n}|\tilde w_n(g)w'(g)|\ \geq\ \de/2
$$
Since $\,H_n\cap G_n=\varnothing\,$, {\bf(a)} implies that
$\,\Ga(H_n)\cap\Ga(H_{n+j})=\varnothing\,$ for all sufficiently
large $\,j\in\N\,$. Therefore we can choose a subsequence
$\,\{H_{n_i}\}_{i=1,2,\ldots}\,$ of the sequence
$\,\{H_n\}_{n=1,2,\ldots}$ in such a way that
\begin{equation}\label{3.1}
\Ga(H_{n_i})\cap\Ga(H_{n_j})=\varnothing\,,\qquad\forall n_i\ne
n_j\,.
\end{equation}

Let $\tilde w(g):=0\,$ for all $\,g\not\in\cup_{i=1}^\infty
H_{n_i}\,$ and $\tilde w(g):=w_{n_i}(g)\,$ for all $\,g\in
H_{n_i}\,$. We have shown in the first part of the proof that the
function $\tilde w\,$ is estimated on every set $\,H_{n_i}\,$ by the
sum of two functions from $\,\PC^0(\Ga)\,$. From here and
\eqref{3.1} it follows that $\tilde w\,$ is estimated by the sum of
two functions from $\,\PC^0(\Ga)\,$ on the whole set $\,\Om\,$.
Therefore $\tilde w\in\WC\,$. On the other hand,
$\,\sum_{g\in\Om}|\tilde w(g)\,w'(g)|=\infty\,$ which contradicts to
the condition $\,w'\in\WC'\,$. This proves (3).
\end{proof}

\begin{remark}\label{R3.2}
From our definition of the operator $\,T_n\,$ it is clear that
$\,T_nw\in\PC(\Ga)\,$ whenever $\,w\,$ takes only the values 0 and
1. This observation can be used for constructing functions
$\,w\in\PC(\Ga)\,$.
\end{remark}

\subsection{Closed convex hull}\label{S3.3}

Lemma \ref{L3.1} allows us to prove the following

\begin{theorem}\label{T3.3}
If the conditions {\bf(a)}, {\bf(w)} and {\bf(a$_2$)} are fulfilled
then
\begin{equation}\label{3.2}
\SC(\Ga,\WC)\ =\ \cconv\ex\SC(\Ga,\WC)
\end{equation}
in any topology $\,\Tf\,$ satisfying the conditions {\bf(w$_1$)} and
{\bf(w$_2$)}.
\end{theorem}

\begin{proof}
Let us fix an arbitrary function $\,w\in\SC(\Ga,\WC)\,$. In view of
{\bf(w$_2$)}, it is sufficient to prove that
$\,w\in\cconv\ex\SC(\Ga,\WC)\,$. Recall that, by the separation
theorem, $\,\cconv\ex\SC(\Ga,\WC)\,$ coincides with the weak closure
of the convex set $\,\conv\ex\SC(\Ga,\WC)\,$.

\smallskip\noindent 1. \
Assume that there exist $\,g_1,g_2,\ldots\in\supp w\,$ such that
$\,\Ga(g_1)=\Ga(g_j)\,$ for all $\,j\geq2\,$. Let
\begin{enumerate}
\item[$\bullet$]
$\,w_i(g_i):=\sum_{j\geq1}w(g_j)\,$,
\item[$\bullet$]
$\,w_i(g_j):=0\,$ for all $\,j\ne i\,$,
\item[$\bullet$]
$\,w_i(g)=w(g)\,$ for all $\,g\not\in\{g_1,g_2,\ldots\}\,$.
\end{enumerate}
Then $\,w_i\in\SC(\Ga,\WC)\,$ and, in view of {\bf(w$_1$)}, there
exists a sequence of finite convex linear combinations of the
functions $\,w_i\,$ which is weakly convergent to $\,w\,$. Since the
set of all finite intersections of the sets $\,\Om_k\,$ is
countable, this implies that $\,w\,$ is contained in the weak
sequential closure of the set of all functions $\,\tilde
w\in\SC(\Ga,\WC)\,$ whose supports satisfy the condition
{\bf(a$_1$)}. Therefore we shall be assuming without loss of
generality that $\,\supp w\,$ satisfies {\bf(a$_1$)}.

\smallskip\noindent 2. \
If $\,\#\,\Ga<\infty\,$ then, by {\bf(a$_1$)}, $\,\#\,\supp
w<\infty\,$. In this case $\,w\,$ belongs to
$\,\conv\ex\SC(\Ga,\WC)\,$. If $\,\#\,\Ga=\infty\,$ and
{\bf(a$_{21}$)} holds then
$\,\SC(\Ga)=\cconv\ex\SC(\Ga)=\varnothing\,$ (see Corollary
\ref{C1.5}). Further on we shall be assuming that
$\,\#\,\Ga=\infty\,$ and $\,\Ga\,$ satisfies {\bf(a$_{22}$)}.

\smallskip\noindent 3. \
Assume that $\,w^\star_n\in\SC_n(\Ga)\,$ and $\,\supp w^\star_n\,$
satisfies the condition {\bf(a$_1$)}. Let $\,\Om_k\cap\supp
w^\star_n=\{g_1^k,g_2^k,...\}\,$ and
$\,\Ga_{i,k}^{(n)}:=\Ga(g_i^k)\cap\{\Om_{n+1},\Om_{n+2},\ldots\}\,$
where $\,k=1,2,\ldots,n\,$. From {\bf(a$_1$)}, {\bf(a)} and
\eqref{1.5} it follows that
$$
\sup_{\Om_j\in\Ga_{i,k}^{(n)}}\;\sum_{g\in\Om_j\cap
G_n}w^\star_n(g)\ \underset{i\to\infty}\to\ 0
$$
for all $\,k=1,2,\ldots,n\,$. Also, if $\,\#\,\Om_k=\infty\,$ then
$$
v_n(g_i^k)\ :=\ \sum_{j>i}w^\star_n(g_j^k)\
\underset{i\to\infty}\to\ 0
$$
and $\,\Ga(g_i^k)\cap\{\Om_1,\Om_2,\ldots,\Om_n\}=\{\Om_k\}\,$
whenever $\,i\,$ is sufficiently large. Therefore for all
sufficiently large $\,i\in\N\,$ the functions $\,w^\star_{n,i}\,$
defined by the equalities
\begin{enumerate}
\item[$\bullet$]
$\,w^\star_{n,i}(g_j^k):=w^\star_n(g_j^k)\,$ for all
$\,j=1,\ldots,i-1\,$,
\item[$\bullet$]
$\,w^\star_{n,i}(g_i^k):=v_n(g_i^k)\,$,
\item[$\bullet$]
$\,w^\star_{n,i}(g_j^k)=0\,$ for all $\,j>i\,$,
\end{enumerate}
belong to $\,\SC_n(\Ga)\,$. Each of these functions has a compact
support and therefore coincides with a finite convex linear
combination of functions $\,w^j_{n,i}\,$ from
$\,\conv\ex\SC_n(\Ga)\,$. By {\bf(w$_1$)}, the sequence
$\,\{T_0w^\star_{n,i}\}_{i=1,2,\ldots}\,$ weakly converges to
$\,T_0w^\star_n\,$. Thus $\,T_0w^\star_n\,$ is contained in the weak
sequential closure of the set $\,T_0(\conv\ex\SC_n(\Ga))\,$.

\smallskip\noindent 4. \
Let $\,w_n^\star:=\left.w\right|_{G_n}\,$. By {\bf(w$_1$)}, we have
$\,\langle w-T_0w_n^\star,w'\rangle\underset{n\to\infty}\to0\,$ for
all $\,w'\in\WC'\,$.

Since $\,\supp w\,$ satisfies {\bf(a$_1$)}, the same is true for
$\,\supp w_n^\star\,$. Let
$\,w_{n,i}^\star:=\sum_j\al_j\,w^j_{n,i}\,$ be the finite convex
linear combinations of functions $\,w^j_{n,i}\in\ex\SC_n(\Ga)\,$
introduced in the previous part of the proof. Then $\,\langle
T_0w_n^\star-T_0w_{n,i},w'\rangle\underset{i\to\infty}\to0\,$ for
each function $\,w'\in\WC'\,$ and for each $\,n\in\N\,$.

Define $\,w_{n,i}:=\sum_j\al_j\,T_nw^j_{n,i}\,$. The conditions (2)
and (3) of Lemma \ref{L3.1} imply that
$\,w_{n,i}\in\conv\ex\SC(\Ga,\WC)\,$ and
\begin{multline*}
\langle T_0w_{n,i}^\star-w_{n,i},w'\rangle\ =\
\sum_j\al_j\,\langle T_0w_{n,i}^j-T_nw_{n,i}^j,w'\rangle\\
\leq\ \sup\limits_{w\in\SC_n^0(\Ga,\WC)} \langle
T_0w-T_nw,w'\rangle\ \underset{n\to\infty}\to\ 0\,,\qquad\forall
w'\in\WC'\,.
\end{multline*}

From the above, it follows that for each function $\,w'\in\WC'\,$,
choosing a sufficiently large $\,n\,$ and then a sufficiently large
$\,i\,$, we can make the sum on the right hand side of the
inequality
\begin{multline*}
|\langle w-w_{n,i},w'\rangle|\leq|\langle
w-T_0w_n^\star,w'\rangle|+|\langle
T_0w_n^\star-T_0w_{n,i}^\star,w'\rangle|+|\langle
T_0w_{n,i}^\star-w_{n,i},w'\rangle|\\
\leq\ |\langle w-T_0w_n^\star,w'\rangle|+|\langle
T_0w_n^\star-T_0w_{n,i}^\star,w'\rangle|+\sup_{w\in\SC_n^0(\Ga,\WC)}
\langle T_0w-T_nw,w'\rangle
\end{multline*}
arbitrarily small. By the separation theorem, this implies
\eqref{3.2}.
\end{proof}

\begin{remark}\label{R3.4}
If we increase the dual space $\,\WC^*\,$ then the topology on
$\,\WC\,$ becomes finer. Therefore, without loss of generality, in
Theorem \ref{T3.3} one can replace {\bf(w$_1$)} with the stronger
condition
\begin{enumerate}
\item[{\bf(w$_1'$)}]
$\,\WC^*=\WC'\,$.
\end{enumerate}
One can also assume that $\,\Tf\,$ is the strongest topology
satisfying {\bf(w$_1'$)} (which is called the Mackey topology).
Finally, since $\,\SC(\Ga)\subset\WC_1\,$, we can always assume that
$\,\WC\subseteq\WC_1\,$ because a reduction of $\,\WC\,$ increases
the space $\,\WC'\,$. If we take a smaller space
$\,\WC\subset\WC_1\,$ (for example, one may wish to consider
functions whose restrictions to $\,\Om_k\,$ belong to a weighted
space $\,l^p\,$) then Theorem \ref{T3.3} gives a stronger result
which is valid for a narrower class of functions
$\,w\in\SC(\Ga,\WC)\,$.
\end{remark}

\begin{remark}\label{R3.5}
The topological space $\,(\WC,\Tf)\,$ may be incomplete. However, if
$\,\WC=\WC''\,$ and $\,\WC^*=\WC'\,$ then $\,\WC\,$ is complete in
the Mackey topology and is weakly sequentially complete (see, for
example, Section 30.5 in \cite{K}).
\end{remark}

\begin{example}\label{E3.6}
Let us consider the coarsest topology on $\,\WC=\WC_1\,$ satisfying
{\bf(w$_1$)} and {\bf(w$_2$)}, with respect to which all the
functionals $\,w\to w(g)\,$ are continuous. This topology is
generated by the seminorms $\,p_k(w):=\sum_{g\in\Om_k}|w(g)|\,$ and
$\,p_g(w):=|w(g)|\,$ and, consequently, is metrizable. Therefore
Theorem \ref{T3.3} implies that, under the conditions {\bf(a)} and
{\bf(a$_2$)}, for every function $\,w\in\SC(\Ga)\,$ there exists a
sequence of functions $\,w_n\in\conv\ex\SC(\Ga)\,$ such that
$$
|w(g)-w_n(g)|\underset{n\to\infty}\to0\quad\text{and}
\quad\sum_{g\in\Om_j}|w(g)-w_n(g)|\underset{n\to\infty}\to0
$$
for all $\,g\in\Om\,$ and $\,k=1,2,\ldots\;$ In \cite{Ke} the
equality \eqref{3.2} was proved for this coarsest topology on the
space of infinite matrices.
\end{example}

A topology $\,\Tf\,$ satisfying the conditions {\bf(w$_1$)} and
{\bf(w$_2$)} (in particular, the Mackey topology), may well be
non-metrizable. Therefore, in the general case, \eqref{3.2} does not
imply the existence of a sequence of convex linear combinations
$\,w_n\in\conv\ex\SC(\Ga,\WC)\,$ convergent to a given function
$\,w\in\SC(\Ga,\WC)\,$. It is possible that Theorem \ref{T3.3} can
be improved in this direction (note that in the parts 1 and 3 of the
proof we spoke about sequential closures).

\begin{remark}\label{R3.7}
It seems to be natural to consider the closure of
$\,\conv\ex\SC(\Ga)\,$ with respect to the norm
$\,\|w\|_\SC:=\sup_j\sum_{g\in\Om_j}|w(g)|\,$ on the space
$\,\WC_1\,$. However, this closure does not always coincide with
$\,\SC(\Ga)\,$.

Indeed, let $\,\Ga\,$ be an infinite family of mutually disjoint
sets $\,\Om_k\,$ and $\,\#\,\Om_k=k\,$. Then
$\,\ex\SC(\Ga)=\PC(\Ga)\,$ and $\,\SC(\Ga)\,$ contains the function
$\,w_0\,$ defined by the equalities
$\,\left.w_0\right|_{\Om_k}\equiv k^{-1}\,$. On the other hand,
\begin{equation}\label{3.3}
\sup_j\#\,\{g\in\Om_j\,:\,w(g)\ne0\}<\infty\,,\qquad\forall
w\in\conv\PC(\Ga)\,.
\end{equation}
Therefore $\,\|w_0-w\|_\SC=1\,$ for all $\,w\in\conv\PC(\Ga)\,$.
\end{remark}


\begin{thebibliography}{MMMMM}

\bibitem[An]{An}
T. Ando. Majorization, doubly stochastic matrices, and comparison
of eigenvalues, {\it Linear Algebra Appl.} {\bf 118} (1989),
163--248.

\bibitem[Bi1]{Bi1}
G. Birkhoff. Three observations on linear algebra. (Spanish) {\it
Univ. Nac. Tucum\'an. Revista A.} {\bf 5} (1946), 147--151.

\bibitem[Bi2]{Bi2}
G. Birkhoff. {\it Lattice Theory.} Amer. Math. Soc. Colloq. Publ.,
vol. 25, rev. ed., Amer. Math. Soc., New York, 1948.

\bibitem[BR]{BR}
R.B. Bapat and T.E.S. Raghavan. {\it Nonnegative matrices and
applications.} Encyclopedia of Mathematics and its Applications,
{\bf 64}. {\sl Cambridge University Press, Cambridge,} 1997.

\bibitem[BS]{BS}
J.R. Brown and R.C. Shiflett. On extreme doubly stochastic
measures. {\it Michigan Math. J.} {\bf 17} (1970), 249--254.

\bibitem[CLMST]{CLMST}
R.M. Caron, X. Li, P. Mikusi\'nski, H.Sherwood and M.D. Taylor.
Nonsquare ``doubly stochastic'' matrices. {\it Distributions with
fixed marginals and related topics (Seattle, WA, 1993),} 187--197,
IMS Lecture Notes Monogr. Ser., {\bf 28}, {\sl Inst. Math.
Statist., Hayward, CA,} 1996.

\bibitem[Do]{Do}
R.G. Douglas. On extreme measures and subspace density. {\it
Michigan Math. J.} {\bf 11} (1964), 243--246.

\bibitem[Fe]{Fe}
D. Feldman. Doubly stochastic measures: three vignettes. {\it
Distributions with fixed marginals and related topics (Seattle,
WA, 1993),} 187--197, IMS Lecture Notes Monogr. Ser., {\bf 28},
{\sl Inst. Math. Statist., Hayward, CA,} 1996.

\bibitem[Gr]{Gr}
R. Grza\'slewicz. On extreme infinite doubly stochastic matrices.
{\it Illinois J. Math.} {\bf 31}, no. 4 (1987), 529--543.

\bibitem[Ho]{Ho}
A.J. Hoffman. A special class of doubly stochastic matrices. {\it
Aequationes Math.} {\bf 2} (1969), 319--326.

\bibitem[Is]{Is}
J.R. Isbell. Birkhoff's problem $111$. {\it Proc. Amer. Math.
Soc.} {\bf 6} (1955), 217--218.

\bibitem[K]{K}
G. K\"othe. {\it Topological Vector Spaces \rm I},
Springer-Verlag, 1969.

\bibitem[Ka]{Ka}
M. Katz. On the extreme points of certain convex polytope. {\it J.
Combinatorial Theory} {\bf 8} (1970), 417--423.

\bibitem[Ke]{Ke}
D.G. Kendall. On infinite doubly-stochastic matrices and
Birkhoff's problem. {\it  J. London Math. Soc.} {\bf 35} (1960),
81--84.

\bibitem[KST]{KST}
A. Kami\'nski, H. Sherwood and M.D. Taylor. Doubly stochastic
measures with mass on the graph of two functions, {\it Real Anal.
Exchange} {\bf 28} (1987/88), no 1, 253--257.

\bibitem[Le]{Le}
G. Letac. Repres\'entation des mesures de probabilit\'e sur le
produit de deux espaces d\'enombrables, de marqes donn\'ees.
(French) {\it Illinois J. Math.} {\bf 10} (1966), 497--507.

\bibitem[Li]{Li}
J. Lindenstrauss. A remark on extreme doubly stochastic measures.
{\it Amer. Math. Monthly} {\bf 72} (1965), 379--382.

\bibitem[LLL]{LLL}
J.L. Lewandowski, C.L. Liu and J.W.-S. Liu. An algorithmic proof
of a generalization of the Birkhoff--von Neumann theorem. {\it J.
Algorithms} {\bf 7} (1986), no. 3, 323--330.

\bibitem[LMST]{LMST}
X. Li, P. Mikusi\'nski, H.Sherwood and M.D. Taylor. In quest of
Birkhoff's theorem in higher dimensions. {\it Distributions with
fixed marginals and related topics (Seattle, WA, 1993),} 187--197,
IMS Lecture Notes Monogr. Ser., {\bf 28}, {\sl Inst. Math.
Statist., Hayward, CA,} 1996.

\bibitem[Mu]{Mu}
H.G. Mukerjee. Supports of extremal measures with given marginals.
(French) {\it Illinois J. Math.} {\bf 29}, no. 2 (1985), 248--260.

\bibitem[Ro]{Ro}
J.V. Romanovsky. A simple proof of the Birkhoff--von Neumann
theorem on bistochastic matrices, {\it A tribute to Ilya Bakelman
(College station, TX, 1993),} 51--53, Discourses math. Appl. {\bf
3}, {\sl Texas A\&M Univ., College Station, TX,} 1994.

\bibitem[RP]{RP}
B.A. Rattray and J.E.L. Peck. Infinite stochastic matrices. {\it
Trans. Roy. Soc. Canada Sect. III (3)} {\bf 49} (1955), 55--57.

\bibitem[Ru]{Ru}
W.H. Ruckle, {\it Sequence Spaces}. Research Notes in Mathematics
{\bf 35}, Pitman Advanced Publishing Program, 1981.

\bibitem[Sa]{Sa}
Yu. Safarov. Birkhoff's theorem and multidimensional spectra of
self-adjoint operators. {\it J. Funct. Anal.} {\bf 222} (2005),
61--97.

\bibitem[ST]{ST}
H. Schreck and G. Tinhofer. A note on certain subpolytopes of the
assignment polytope associated with circulant graphs. {\it Linear
Algebra Appl.} {\bf 111} (1988), 125--134.

\bibitem[Ti]{Ti}
G. Tinhofer. Strong tree-cographs are Birkhoff graphs. {\it
Discrete Appl. Math} {\bf 22} (1988/89), no. 3, 275--288.

\bibitem[Vi]{Vi}
R.A. Vitale. Parametrizing doubly stochastic measures. {\it
Distributions with fixed marginals and related topics (Seattle,
WA, 1993),} 187--197, IMS Lecture Notes Monogr. Ser., {\bf 28},
{\sl Inst. Math. Statist., Hayward, CA,} 1996.

\end{thebibliography}
\end{document}